\documentclass[reqno]{amsart}%
\usepackage{amsmath}
\usepackage{amsfonts}
\usepackage{amssymb}
\usepackage{graphicx}%
\newtheorem{theorem}{Theorem}

\newtheorem{lemma}[theorem]{Lemma}

\begin{document}
\title[Symmetry of Solutions]{Symmetry of Solutions for Some Bernoulli Initial-Value Problems}
\author{Nadejda E. Dyakevich}
\address{Department of Mathematics, California State University San Bernardino, 5500
University Parkway, San Bernardino, CA, 92407-2397, USA}
\email{dyakevic@csusb.edu}
\keywords{Bernoulli equation, symmetry, graphs.}
\subjclass{34A34, 34C14}
\maketitle

\begin{abstract}
In this paper, we show that under certain conditions on the coefficients and
initial values, solutions of two different Bernoulli initial-value problems
are symmetric to each other either with respect to the $t$-axis, or the
$y$-axis, or the origin. This knowledge of symmetry is very useful since
solving one Bernoulli problem solves automatically the other. To the best of
our knowledge, this interesting property has not been described in any other papers.

\end{abstract}

\section{ Introduction}

Let us consider the following two Bernoulli initial-value problems:%
\begin{equation}
\left\{
\begin{array}
[c]{c}%
\dfrac{dy_{1}(t)}{dt}=a_{1}(t)y_{1}(t)+b_{1}(t)y_{1}^{n}(t)\\
y_{1}(0)=d_{1}\text{,}%
\end{array}
\right.  \tag{1}\label{1}%
\end{equation}
and
\begin{equation}
\left\{
\begin{array}
[c]{c}%
\dfrac{dy_{2}(t)}{dt}=a_{2}(t)y_{2}(t)+b_{2}(t)y_{2}^{n}(t)\\
y_{2}(0)=d_{2}\text{,}%
\end{array}
\right.  \tag{2}\label{2}%
\end{equation}
where $d_{1}$ and $d_{2}$ are nonzero real numbers, $n$ is a rational number,
and the functions $a_{1}(t)$, $b_{1}(t)$, $a_{2}(t)$, and $b_{2}(t)$ are
continuous on some interval that contains $t=0$. Then the problems (\ref{1})
and (\ref{2}) have unique solutions $y_{1}(t)$ and $y_{2}(t)$ in some
intervals $I_{1}$ and $I_{2}$, respectively ([1] - [6]). The results in this
paper are also valid for linear ($n=0$ and $n=1$) and for Riccati ($n=2$)
initial-value problems.

If $n=\frac{p}{q}$, where $p$ is an even number and $q$ is an odd number, then
the solutions are:%
\begin{equation}
y_{i}(t)=e^{%
{\displaystyle\int\nolimits_{0}^{t}}
a_{i}(\hat{t})d\hat{t}}\left(  \frac{1}{d_{i}^{n-1}}-\left(  n-1\right)
{\displaystyle\int\nolimits_{0}^{t}}
b_{i}(\tilde{t})e^{\left(  n-1\right)
{\displaystyle\int\nolimits_{0}^{\tilde{t}}}
a_{i}(\hat{t})d\hat{t}}d\tilde{t}\right)  ^{-\frac{1}{n-1}}\text{,}
\tag{3}\label{3}%
\end{equation}
where $i=1$, $2$. If $n=\frac{p}{q}$, where $p$ is an odd number and $q$ is an
even number, then (\ref{3}) is valid if the initial values $d_{1}$ and $d_{2}$
are positive. If $n=\frac{p}{q}\neq1$, where $p$ and $q$ are odd numbers, then%
\[
y_{i}(t)=\pm e^{%
{\displaystyle\int\nolimits_{0}^{t}}
a_{i}(\hat{t})d\hat{t}}\left(  \frac{1}{d_{i}^{n-1}}-\left(  n-1\right)
{\displaystyle\int\nolimits_{0}^{t}}
b_{i}(\tilde{t})e^{\left(  n-1\right)
{\displaystyle\int\nolimits_{0}^{\tilde{t}}}
a_{i}(\hat{t})d\hat{t}}d\tilde{t}\right)  ^{-\frac{1}{n-1}}\text{,}%
\]
where $i=1$, $2$ and the sign is the same as the sign of the corresponding
initial value $d_{i}$, $i=1$, $2$. If $n=1$, then we have:%
\[
y_{i}(t)=d_{i}e^{%
{\displaystyle\int\nolimits_{0}^{t}}
\left(  a_{i}(\hat{t})+b_{i}(\hat{t})\right)  d\hat{t}}\text{,}%
\]
where $i=1$, $2$.

The function $y_{1}(t)$ is symmetric to $y_{2}(t)$ with respect to the origin,
if $y_{2}(-t)=-y_{1}(t)$. Also, $y_{1}(t)$ is symmetric to $y_{2}(t)$ with
respect to the $t$-axis, if $y_{2}(t)=-y_{1}(t)$. Finally, $y_{1}(t)$ is
symmetric to $y_{2}(t)$ with respect to the $y$-axis, if $y_{2}(-t)=y_{1}(t)$.

In Section 2, we show that under certain conditions on the
coefficients $a_{1}(t)$, $b_{1}(t)$, $a_{2}(t)$, and $b_{2}(t)$,
and the initial values $d_{1}$ and $d_{2}$, the solutions
$y_{1}(t)$ and $y_{2}(t)$ are symmetric to each other either with
respect to the $t$-axis, or the $y$-axis, or the origin. This
knowledge of symmetry can be very useful since solving one
Bernoulli problem solves automatically the other. For example, the
results in this paper can be used to prepare multiple versions of
examination questions. In addition, solutions for Bernoulli
problems are used to study properties and long term behavior of
solutions for more complicated ordinary and partial differential
equations.

\section{Symmetry of solutions}

\begin{lemma}
Let $n$\ be any real number, $a(t)$\ be an even function, and $b(t)$\ be an
odd function. Then the following holds:%
\begin{equation}%
{\displaystyle\int\nolimits_{-t}^{0}}
b(\tilde{t})e^{-\left(  n-1\right)
{\displaystyle\int\nolimits_{0}^{\tilde{t}}}
a(\hat{t})d\hat{t}}d\tilde{t}=-%
{\displaystyle\int\nolimits_{0}^{t}}
b(\tilde{t})e^{\left(  n-1\right)
{\displaystyle\int\nolimits_{0}^{\tilde{t}}}
a(\hat{t})d\hat{t}}d\tilde{t}\text{.} \tag{4}\label{4}%
\end{equation}

\begin{proof}
If $a(t)$ is an even function, then
\[
\left(  n-1\right)
{\displaystyle\int\nolimits_{0}^{\tilde{t}}}
a(\hat{t})d\hat{t}%
\]
is an odd function. Therefore,%
\[
e^{\left(  n-1\right)
{\displaystyle\int\nolimits_{0}^{-\tilde{t}}}
a(\hat{t})d\hat{t}}=e^{-\left(  n-1\right)
{\displaystyle\int\nolimits_{0}^{\tilde{t}}}
a(\hat{t})d\hat{t}}\text{.}%
\]
The graphs of the functions
\[
f(\tilde{t})=e^{-\left(  n-1\right)
{\displaystyle\int\nolimits_{0}^{\tilde{t}}}
a(\hat{t})d\hat{t}}%
\]
and
\[
g(\tilde{t})=e^{\left(  n-1\right)
{\displaystyle\int\nolimits_{0}^{\tilde{t}}}
a(\hat{t})d\hat{t}}%
\]
are reflections of each other with respect to the $y$-axis since%
\begin{align}
g(-\tilde{t})  &  =e^{\left(  n-1\right)
{\displaystyle\int\nolimits_{0}^{-\tilde{t}}}
a(\hat{t})d\hat{t}}=e^{-\left(  n-1\right)
{\displaystyle\int\nolimits_{-\tilde{t}}^{0}}
a(\hat{t})d\hat{t}}\tag{5}\label{5}\\
&  =e^{-\left(  n-1\right)
{\displaystyle\int\nolimits_{0}^{\tilde{t}}}
a(\hat{t})d\hat{t}}=f(\tilde{t})\text{.}\nonumber
\end{align}
Let us consider the functions $h(\tilde{t})=b(\tilde{t})g(\tilde{t})$ and
$s(\tilde{t})=b(\tilde{t})f(\tilde{t})$. Using (\ref{5}), we have:%
\begin{equation}
h(-\tilde{t})=b(-\tilde{t})g(-\tilde{t})=-b(\tilde{t})f(\tilde{t}%
)=-s(\tilde{t})\text{,} \tag{6}\label{6}%
\end{equation}
which implies (\ref{4}).
\end{proof}
\end{lemma}

\begin{theorem}
Let $n$\ $=\frac{p}{q}$, where $p$\ is an even number and $q$\ is an odd
number. The following holds for the solutions $y_{1}(t)$\ and $y_{2}(t)$\ of
the problems (\ref{1}) and (\ref{2}):

(i). Let $a_{1}(t)$\ be an even function, $b_{1}(t)$\ be an odd function, and
$a_{2}(t)=-a_{1}(t)$, $b_{2}(t)=-b_{1}(t)$, and $d_{2}=-d_{1}$. Then
$y_{1}(t)$\ and $y_{2}(t)$\ are symmetric to each other with respect to the origin.

(ii).Let $a_{1}(t)$\ be an odd function, $b_{1}(t)$\ be an even function, and
$a_{2}(t)=a_{1}(t)$, $b_{2}(t)=b_{1}(t)$, and $d_{2}=-d_{1}$. Then $y_{1}%
(t)$\ and $y_{2}(t)$\ are symmetric to each other with respect to the origin.

(iii).Let $a_{1}(t)$\ and $b_{1}(t)$\ be even functions, and $a_{2}%
(t)=-a_{1}(t)$, $b_{2}(t)=b_{1}(t)$, and $d_{2}=-d_{1}$. Then $y_{1}(t)$\ and
$y_{2}(t)$\ are symmetric to each other with respect to the origin.

(iv). If $a_{1}(t)=a_{2}(t)$, $b_{2}(t)=-b_{1}(t)$, and $d_{2}=-d_{1}$, then
$y_{1}(t)$\ and $y_{2}(t)$\ are symmetric to each other with respect to the
$t$-axis.

\begin{proof}
(i). Let us show that $y_{2}(-t)=-y_{1}(t)$. Using (\ref{4}), we have:%
\begin{align*}
y_{2}(-t)  &  =e^{-%
{\displaystyle\int\nolimits_{0}^{-t}}
a_{1}(\hat{t})d\hat{t}}\left(  -\frac{1}{d_{1}^{n-1}}+\left(  n-1\right)
{\displaystyle\int\nolimits_{0}^{-t}}
b_{1}(\tilde{t})e^{-\left(  n-1\right)
{\displaystyle\int\nolimits_{0}^{\tilde{t}}}
a_{1}(\hat{t})d\hat{t}}d\tilde{t}\right)  ^{-\frac{1}{n-1}}\\
&  =-y_{1}(t)\text{.}%
\end{align*}
(ii). If $a_{1}(t)$ is an odd function, and $b_{1}(t)$ is an even function,
then
\[
b_{1}(t)e^{\left(  n-1\right)
{\displaystyle\int\nolimits_{0}^{t}}
a_{1}(\hat{t})d\hat{t}}%
\]
is even. This implies that%
\begin{equation}%
{\displaystyle\int\nolimits_{-t}^{0}}
b_{1}(\tilde{t})e^{\left(  n-1\right)
{\displaystyle\int\nolimits_{0}^{\tilde{t}}}
a_{1}(\hat{t})d\hat{t}}d\tilde{t}=%
{\displaystyle\int\nolimits_{0}^{t}}
b_{1}(\tilde{t})e^{\left(  n-1\right)
{\displaystyle\int\nolimits_{0}^{\tilde{t}}}
a_{1}(\hat{t})d\hat{t}}d\tilde{t}\text{.} \tag{7}\label{7}%
\end{equation}
Using (\ref{7}), we have:%
\begin{align*}
y_{2}(-t)  &  =e^{%
{\displaystyle\int\nolimits_{0}^{-t}}
a_{1}(\hat{t})d\hat{t}}\left(  -\frac{1}{d_{1}^{n-1}}-\left(  n-1\right)
{\displaystyle\int\nolimits_{0}^{-t}}
b_{1}(\tilde{t})e^{\left(  n-1\right)
{\displaystyle\int\nolimits_{0}^{\tilde{t}}}
a_{1}(\hat{t})d\hat{t}}d\tilde{t}\right)  ^{-\frac{1}{n-1}}\\
&  =-y_{1}(t)\text{.}%
\end{align*}
(iii). The functions
\[
f(t)=b_{1}(t)e^{\left(  n-1\right)
{\displaystyle\int\nolimits_{0}^{t}}
a_{1}(\hat{t})d\hat{t}}%
\]
and
\[
g(t)=b_{1}(t)e^{-\left(  n-1\right)
{\displaystyle\int\nolimits_{0}^{t}}
a_{1}(\hat{t})d\hat{t}}%
\]
are symmetric to each other with respect to the $y$-axis because%
\[
f(-t)=b_{1}(-t)e^{\left(  n-1\right)
{\displaystyle\int\nolimits_{0}^{-t}}
a_{1}(\hat{t})d\hat{t}}=b_{1}(t)e^{-\left(  n-1\right)
{\displaystyle\int\nolimits_{0}^{t}}
a_{1}(\hat{t})d\hat{t}}=g(t)\text{.}%
\]
Therefore,%
\begin{equation}%
{\displaystyle\int\nolimits_{-t}^{0}}
b_{1}(\tilde{t})e^{\left(  n-1\right)
{\displaystyle\int\nolimits_{0}^{\tilde{t}}}
a_{1}(\hat{t})d\hat{t}}d\tilde{t}=%
{\displaystyle\int\nolimits_{0}^{t}}
b_{1}(\tilde{t})e^{-\left(  n-1\right)
{\displaystyle\int\nolimits_{0}^{\tilde{t}}}
a_{1}(\hat{t})d\hat{t}}d\tilde{t}\text{.} \tag{8}\label{8}%
\end{equation}
and%
\begin{equation}%
{\displaystyle\int\nolimits_{-t}^{0}}
b_{1}(\tilde{t})e^{-\left(  n-1\right)
{\displaystyle\int\nolimits_{0}^{\tilde{t}}}
a_{1}(\hat{t})d\hat{t}}d\tilde{t}=%
{\displaystyle\int\nolimits_{0}^{t}}
b_{1}(\tilde{t})e^{\left(  n-1\right)
{\displaystyle\int\nolimits_{0}^{\tilde{t}}}
a_{1}(\hat{t})d\hat{t}}d\tilde{t}\text{.} \tag{9}\label{9}%
\end{equation}
Using (\ref{9}), we have:%
\begin{align*}
y_{2}(-t)  &  =e^{-%
{\displaystyle\int\nolimits_{0}^{-t}}
a_{1}(\hat{t})d\hat{t}}\left(  -\frac{1}{d_{1}^{n-1}}-\left(  n-1\right)
{\displaystyle\int\nolimits_{0}^{-t}}
b_{1}(\tilde{t})e^{-\left(  n-1\right)
{\displaystyle\int\nolimits_{0}^{\tilde{t}}}
a_{1}(\hat{t})d\hat{t}}d\tilde{t}\right)  ^{-\frac{1}{n-1}}\\
&  =-y_{1}(t)\text{.}%
\end{align*}
(iv). Using (\ref{3}), we have:%
\begin{align*}
y_{2}(t)  &  =e^{%
{\displaystyle\int\nolimits_{0}^{t}}
a_{1}(\hat{t})d\hat{t}}\left(  -\frac{1}{d_{1}^{n-1}}+\left(  n-1\right)
{\displaystyle\int\nolimits_{0}^{t}}
b_{1}(\tilde{t})e^{\left(  n-1\right)
{\displaystyle\int\nolimits_{0}^{\tilde{t}}}
a_{1}(\hat{t})d\hat{t}}d\tilde{t}\right)  ^{-\frac{1}{n-1}}\\
&  =-y_{1}(t)\text{.}%
\end{align*}

\end{proof}
\end{theorem}

\begin{theorem}
Let $n$\ be any rational number. Suppose that one of the following is true:

(i). Let $a_{1}(t)$\ be an even function, and $b_{1}(t)$\ be an odd function.
Also, let $a_{2}(t)=-a_{1}(t)$, $b_{2}(t)=b_{1}(t)$, and $d_{2}=d_{1}$.

(ii). Let $a_{1}(t)$\ be an odd function, and $b_{1}(t)$\ be an even function.
Also, let $a_{2}(t)=a_{1}(t)$, $b_{2}(t)=-b_{1}(t)$, and $d_{2}=d_{1}$.

(iii).Let $a_{1}(t)$\ and $b_{1}(t)$\ be even functions. Also, let
$a_{2}(t)=-a_{1}(t)$, $b_{2}(t)=-b_{1}(t)$, and $d_{2}=d_{1}$.

Then the solutions $y_{1}(t)$\ and $y_{2}(t)$\ of the problems (\ref{1}) and
(\ref{2}) are symmetric to each other with respect to the $y$-axis.

\begin{proof}
(i). Suppose that $d_{i}>0$, $i=1$, $2$ and $n\neq1$. Using (\ref{4}), we
have:
\begin{align*}
y_{2}(-t)  &  =e^{-%
{\displaystyle\int\nolimits_{0}^{-t}}
a_{1}(\hat{t})d\hat{t}}\left(  \frac{1}{d_{1}^{n-1}}-\left(  n-1\right)
{\displaystyle\int\nolimits_{0}^{-t}}
b_{1}(\tilde{t})e^{-\left(  n-1\right)
{\displaystyle\int\nolimits_{0}^{\tilde{t}}}
a_{1}(\hat{t})d\hat{t}}d\tilde{t}\right)  ^{-\frac{1}{n-1}}\\
&  =y_{1}(t)\text{.}%
\end{align*}
The above result is also true for the case $d_{i}<0$, $i=1$, $2$, provided
$p$\ is even and $q$\ is odd, or $p$\ and $q$\ are both odd. If $n=1$, then we
have:%
\[
y_{2}(-t)=d_{1}e^{-%
{\displaystyle\int\nolimits_{-t}^{0}}
\left(  -a_{1}(\hat{t})+b_{1}(\hat{t})\right)  d\hat{t}}=y_{1}\left(
t\right)  \text{.}%
\]
(ii). The proof is similar to the proof of Theorem 3 (i), with (\ref{7}) used
instead of (\ref{4}).\newline(iii). The proof is similar to the proof in
Theorem 3 (i), with (\ref{8}) used instead of (\ref{4}).\newline
\end{proof}
\end{theorem}

\begin{theorem}
Let $n=\frac{p}{q}$, where $p$\ and $q$\ are odd numbers.\ The following holds
for the solutions $y_{1}(t)$\ and $y_{2}(t)$\ of the problems (\ref{1}) and
(\ref{2}):

(i). Let $a_{1}(t)$\ be an odd function, $b_{1}(t)$\ be an even function, and
$a_{2}(t)=a_{1}(t)$, $b_{2}(t)=-b_{1}(t)$, and $d_{2}=-d_{1}$. Then $y_{1}%
(t)$\ and $y_{2}(t)$\ are symmetric to each other with respect to the origin.

(ii). Let $a_{2}(t)=a_{1}(t)$, $b_{2}(t)=b_{1}(t)$, and $d_{2}=-d_{1}$. Then
$y_{1}(t)$\ and $y_{2}(t)$\ are symmetric to each other with respect to the
$t$-axis.

(iii). Let $a_{1}(t)$\ be an even function, $b_{1}(t)$\ be an odd function,
and $a_{2}(t)=-a_{1}(t)$, $b_{2}(t)=b_{1}(t)$, and $d_{2}=-d_{1}$. Then
$y_{1}(t)$\ and $y_{2}(t)$\ are symmetric to each other with respect to the origin.

(iv). Let $a_{1}(t)$\ and $b_{1}(t)$\ be even functions, and $a_{2}%
(t)=-a_{1}(t)$, $b_{2}(t)=-b_{1}(t)$, and $d_{2}=-d_{1}$. Then $y_{1}(t)$\ and
$y_{2}(t)$\ are symmetric to each other with respect to the origin.

\begin{proof}
Suppose that $d_{1}<0$ and $d_{2}>0$. (The case $d_{1}>0$ and $d_{2}<0$ can be
shown similarly.) For $n\neq1$,%
\[
y_{1}(t)=-e^{%
{\displaystyle\int\nolimits_{0}^{t}}
a_{1}(\hat{t})d\hat{t}}\left(  \frac{1}{d_{1}^{n-1}}-\left(  n-1\right)
{\displaystyle\int\nolimits_{0}^{t}}
b_{1}(\tilde{t})e^{\left(  n-1\right)
{\displaystyle\int\nolimits_{0}^{\tilde{t}}}
a_{1}(\hat{t})d\hat{t}}d\tilde{t}\right)  ^{-\frac{1}{n-1}}\text{.}%
\]
(i). If $a_{1}(t)$ is an odd function, and $b_{1}(t)$ is an even function,
then the function
\[
b_{1}(t)e^{\left(  n-1\right)
{\displaystyle\int\nolimits_{0}^{t}}
a_{1}(\hat{t})d\hat{t}}%
\]
is even. This implies that%
\begin{equation}%
{\displaystyle\int\nolimits_{-t}^{0}}
b_{1}(\tilde{t})e^{\left(  n-1\right)
{\displaystyle\int\nolimits_{0}^{\tilde{t}}}
a_{1}(\hat{t})d\hat{t}}d\tilde{t}=%
{\displaystyle\int\nolimits_{0}^{t}}
b_{1}(\tilde{t})e^{\left(  n-1\right)
{\displaystyle\int\nolimits_{0}^{\tilde{t}}}
a_{1}(\hat{t})d\hat{t}}d\tilde{t}\text{.}\tag{10}\label{10}%
\end{equation}
Using (\ref{10}), we have for $n\neq1$:%
\begin{align*}
y_{2}(-t) &  =e^{%
{\displaystyle\int\nolimits_{0}^{-t}}
a_{1}(\hat{t})d\hat{t}}\left(  \frac{1}{d_{1}^{n-1}}+\left(  n-1\right)
{\displaystyle\int\nolimits_{0}^{-t}}
b_{1}(\tilde{t})e^{\left(  n-1\right)
{\displaystyle\int\nolimits_{0}^{\tilde{t}}}
a_{1}(\hat{t})d\hat{t}}d\tilde{t}\right)  ^{-\frac{1}{n-1}}\\
&  =-y_{1}(t)\text{.}%
\end{align*}
If $n=1$, then:%
\[
y_{2}(-t)=d_{2}e^{%
{\displaystyle\int\nolimits_{0}^{-t}}
\left(  a_{2}(\hat{t})+b_{2}(\hat{t})\right)  d\hat{t}}=-d_{1}e^{%
{\displaystyle\int\nolimits_{0}^{-t}}
\left(  a_{1}(\hat{t})-b_{1}(\hat{t})\right)  d\hat{t}}=-y_{1}\left(
t\right)  \text{.}%
\]
(ii). We have for $n\neq1$:%
\begin{align*}
y_{2}(t) &  =e^{%
{\displaystyle\int\nolimits_{0}^{t}}
a_{1}(\hat{t})d\hat{t}}\left(  \frac{1}{d_{1}^{n-1}}-\left(  n-1\right)
{\displaystyle\int\nolimits_{0}^{t}}
b_{1}(\tilde{t})e^{\left(  n-1\right)
{\displaystyle\int\nolimits_{0}^{\tilde{t}}}
a_{1}(\hat{t})d\hat{t}}d\tilde{t}\right)  ^{-\frac{1}{n-1}}\\
&  =-y_{1}(t)\text{. }%
\end{align*}
If $n=1$, then:%
\[
y_{2}(t)=d_{2}e^{%
{\displaystyle\int\nolimits_{0}^{t}}
\left(  a_{2}(\hat{t})+b_{2}(\hat{t})\right)  d\hat{t}}=-d_{1}e^{%
{\displaystyle\int\nolimits_{0}^{t}}
\left(  a_{1}(\hat{t})+b_{1}(\hat{t})\right)  d\hat{t}}=-y_{1}\left(
t\right)  \text{.}%
\]
(iii). The proof is similar to the proof of Theorem 4 (i), with (\ref{4}) used
instead of (\ref{10}).\newline(iv). The proof is similar to the proof of
Theorem 4 (i), with (\ref{9}) used instead of (\ref{10}).
\end{proof}
\end{theorem}

\end{document}